\documentclass[journal]{IEEEtran}

\usepackage{balance}
\usepackage{float}
\usepackage{caption}
\usepackage{subcaption}
\usepackage{tikz,pgfplots,pgfplotstable,booktabs}
\usepackage[official]{eurosym}
\usepackage{xcolor}
\usepackage{colortbl}
\usepackage{comment}
\usepackage{cite}
\usepackage{circuitikz}
\usepackage{amsmath,amsfonts,amsthm,bm}
\allowdisplaybreaks[1]
\usepgfplotslibrary{fillbetween}
\usepgfplotslibrary{statistics}

\pgfplotsset{
    discard if not/.style 2 args={
        x filter/.code={
            \edef\tempa{\thisrow{#1}}
            \edef\tempb{#2}
            \ifx\tempa\tempb
            \else
                
            \fi
    }},
    boxplot/hide outliers/.code={
        \def\pgfplotsplothandlerboxplot@outlier{}}
    }
        
\newenvironment{ldescription}[1]
  {\begin{list}{}%
   {\renewcommand\makelabel[1]{##1\hfill}%
   \settowidth\labelwidth{\makelabel{#1}}%
   \setlength\leftmargin{\labelwidth}
   \addtolength\leftmargin{\labelsep}}}
  {\end{list}}
  
\newcommand{\pg}{p^{G}_{it}}
\newcommand{\qg}{q^{G}_{it}}
\newcommand{\pd}{p^{D}_{jt}}
\newcommand{\qd}{q^{D}_{jt}}
\newcommand{\pf}{p^{F}_{lt}}
\newcommand{\qf}{q^{F}_{lt}}

\begin{document}

\title{Learning the price response of active distribution networks for TSO-DSO coordination}

\author{J. M. Morales, S. Pineda and Y. Dvorkin

\thanks{J. M. Morales is with the Department of Applied Mathematics, University of M\'alaga, M\'alaga, Spain. E-mail: juan.morales@uma.es}
\thanks{S. Pineda is with the Department
of Electrical Engineering, University of M\'alaga, M\'alaga, Spain. E-mail: spinedamorente@gmail.com.}
\thanks{Y. Dvorkin is with the New York University, Brooklyn, NY 11201 USA. E-mail: dvorkin@nyu.edu.}

\thanks{This work was supported in part by the European
Research Council (ERC) under the EU Horizon 2020 research and innovation
program (grant agreement No. 755705), in part by the Spanish Ministry of Science and Innovation through project PID2020-115460GB-I00, in part by the Andalusian Regional Government through project P20-00153, and in part by the European Regional Development Fund (FEDER) through the research project UMA2018-FEDERJA-150. The authors thankfully acknowledge the computer resources, technical expertise and assistance provided by the SCBI (Supercomputing and Bioinformatics) center of the University of M\'alaga.}}

\maketitle

\begin{abstract}
The increase in distributed energy resources and flexible electricity consumers has turned TSO-DSO coordination strategies into a challenging problem. Existing decomposition/decentralized methods apply divide-and-conquer strategies to trim down the computational burden of this complex problem, but rely on access to proprietary information or fail-safe real-time communication infrastructures. To overcome these drawbacks, we propose in this paper a TSO-DSO coordination strategy that only needs a series of observations of the nodal price and the power intake at the substations connecting the transmission and distribution networks. Using this information, we learn the price response of active distribution networks (DN) using a decreasing step-wise function that can also adapt to some contextual information.
The learning task can be carried out in a computationally efficient manner and the curve it produces can be interpreted as a market bid, thus averting the need to revise the current operational procedures for the transmission network. Inaccuracies derived from the learning task may lead to suboptimal decisions. However, results from a realistic case study show that the proposed methodology yields operating decisions very close to those obtained by a fully centralized coordination of transmission and distribution.
\end{abstract}


\IEEEpeerreviewmaketitle

\section*{Nomenclature}
The main symbols used throughout this paper are listed next for quick reference. Others are defined as required in the text. We consider hourly time steps and, hence, MW and MWh are interchangeable under this premise.


\subsection{Indexes and sets}
\begin{ldescription}{$xxxxx$}
\item[$b$] Index of blocks in the step-wise approximation of the price response function.
\item[$i$] Index of generating units.
\item[$j$] Index of consumers.
\item[$k$] Index of distribution networks.
\item[$l$] Index of lines.
\item[$n$] Index of nodes.
\item[$t$] Index of time periods. 
\item[$\mathcal{B}$] Set of blocks in the step-wise approximation of the price response function.
\item[$I^T$] Set of generating units at the transmission network.
\item[$I_n$] Set of generating units at node $n$.
\item[$I^D_k$] Set of generating units of distribution network $k$.
\item[$J^T$] Set of consumers in the transmission network.
\item[$J_n$] Set of consumers at node $n$.
\item[$J^D_k$] Set of consumers in distribution network $k$.
\item[$K^T$] Set of distribution networks.
\item[$K_n$] Set of distribution networks at node $n$.
\item[$L^T$] Set of lines of the transmission network.
\item[$L^D_k$] Set of lines of distribution network $k$.
\item[$N^T$] Set of nodes of the transmission network.
\item[$N^D_k$] Set of nodes of distribution network $k$.
\item[$\mathcal{T}$] Set of time periods.
\item[$\mathcal{T}^C(t)$] Subset of time periods that are the closest to $t$.
\end{ldescription}
\subsection{Parameters}
\begin{ldescription}{$xxx$}
\item[$a_i$] Quadratic cost parameter of unit $i$ [\euro/MW$^2$].
\item[$b_i$] Linear cost parameter of unit $i$ [\euro/MW].
\item[$\underline{p}^G_{i}$] Minimum active power output of unit $i$ [MW].
\item[$\overline{p}^G_{i}$] Maximum active power output of  unit $i$ [MW].
\item[$\widehat{p}^D_{jt}$] Baseline demand of consumer $j$ at time $t$ [MW].
\item[$\overline{p}^D_{jt}$] Maximum demand of consumer $j$ at time $t$ [MW].
\item[$\underline{p}^D_{jt}$] Minimum demand of consumer $j$ at time $t$ [MW].
\item[$\underline{q}^G_{i}$] Minimum reactive power output of unit $i$ [MW].
\item[$\overline{q}^G_{i}$] Maximum reactive power output of unit $i$ [MW].
\item[$r_l$] Resistance of line $l$ [p.u.].
\item[$\overline{s}^F_l$] Capacity of line $l$ [MVA].
\item[$\overline{s}^G_i$] Inverter power rate of unit $i$ [MVA].
\item[$s^B$] Base power [MVA].
\item[$\overline{u}^B_{ktb}$] Marginal utility of block $b$ for DN $k$ and time $t$ [\euro/MW].
\item[$\overline{v}_n$] Maximum squared voltage at node $n$ [p.u.].
\item[$\underline{v}_n$] Minimum squared voltage at node $n$ [p.u.].
\item[$x_l$] Reactance of line $l$ [p.u.].
\item[$\alpha_{jt}$] Intercept of the inverse demand function of consumer $j$ at time $t$ [MW].
\item[$\beta_{jt}$] Slope of the inverse demand function of consumer $j$ at time $t$ [\euro/MW$^2$].
\item[$\gamma_j$] Power factor of consumer $j$.
\item[$\delta_j$] Flexibility parameter of consumer $j$.
\item[$\lambda_{kt}$] Price at the substation of DN $k$ at time $t$ [\euro/MW].
\item[$\rho_{it}$] Capacity factor of unit $i$ at time $t$.
\item[$\chi_{kt}$] Contextual information of DN $k$ at time $t$.
\end{ldescription}

\subsection{Variables}
\begin{ldescription}{$xxx$}
\item[$p^G_{it}$] Active power output of unit $i$ at time $t$ [MW].
\item[$p^D_{jt}$] Active power demand of consumer $j$ at time $t$ [MW].
\item[$p^F_{lt}$] Active power flow through line $l$ at time $t$ [MW].
\item[$p^N_{kt}$] Active power intake of DN $k$ at time $t$ [MW].
\item[$q^G_{it}$] Reactive power output of unit $i$ at time $t$ [MVAr].
\item[$q^D_{jt}$] Reactive power demand of consumer $j$ at time $t$ [MVAr].
\item[$q^F_{lt}$] Reactive power flow through line $l$ at time $t$ [MVAr].
\item[$v_{nt}$] Squared voltage magnitude at node $n$ and time $t$ [p.u.].
\item[$\theta_{nt}$] Voltage angle at node $n$ and time $t$ [rad].
\end{ldescription}

\section{Introduction} \label{sec:introduction}

\IEEEPARstart{E}{lectric} power distribution  has been traditionally ignored in the operation of transmission power networks, on the grounds that distribution grids only housed passive loads. However, the proliferation of distributed energy resources (DERs) is rendering this traditional \emph{modus operandi} obsolete \cite{Li2016a}. Power systems engineers are faced with  an unprecedented challenge of efficiently integrating a vast number and a wide spectrum of flexible power assets located in mid- and low-voltage networks into the operation of the transmission power network \cite{Li2018}. Naturally, succeeding in this endeavor  requires the coordination between the transmission and distribution system operators (TSOs and DSOs, respectively), all united in the purpose of fostering an active role of DERs in the operation of the power system through their participation in wholesale electricity markets.

As a result, research emphasis is placed on mechanisms that strengthen the TSO-DSO coordination so that the available flexibility of DERs can be harvested for transmission and wholesale market services \cite{migliavacca2017smartnet,de2019control}. For instance, some recent research works investigate TSO-DSO coordination schemes to improve voltage stability \cite{li2017impact,valverde2019coordination}. Other authors focus on the economic coordination between transmission and distribution operations to minimize total system costs \cite{li2018new}. The present work belongs to this latter group.

Regarding TSO-DSO economic coordination, a single centralized operational model that includes both transmission and distribution networks with their full level of detail is not viable due to its computational cost, modeling complexity and potential conflict of interests between the involved parties \cite{GERARD201840}. Rather, the coordination of transmission and distribution power assets calls for a divide-and-conquer strategy that alleviates the computational burden, allows for decentralization and minimizes the need for information exchange between the TSO and DSOs \cite{kargarian2018toward}. For instance, authors of \cite{yuan_2017} uses Benders decomposition to find the optimal economic dispatch considering TSO-DSO interactions. Similarly, reference \cite{bragin_2017} proposes
a model to operate transmission and distribution systems
in a coordinated manner using a surrogate Lagrangian Relaxation approach. Finally, an analytical target cascading procedure (ATC) to coordinate the operation of transmission and distribution networks is described in \cite{nawaz_2020}.

The decomposition and decentralized methods previously described are able to obtain the same solution as the centralized approach while significantly reducing the computational burden. Yet these methods also have meaningful drawbacks. Decomposition methods still require full access to all physical and economic information on distribution networks. However, as stated in \cite{mohammadi2019diagonal} ``distribution system operators are autonomous entities that are unwilling to reveal their commercially sensitive information.'' Therefore, these methods can hardly be accommodated in a real-life distribution  environment with even a few ambiguous or unknown  parameters (e.g.  topological configuration,  impedance, voltage and flow limits) and proprietary customer-end and  behind-the-meter parameters (e.g. production/utility cost functions, supply/demand elasticity and behavioral aspects of electricity demand). Similarly, decentralized methods are based on repetitive real-time information exchanges between the TSO and the DSOs, and thus, rely on robust and fast communication infrastructure. As discussed in \cite{kargarian2018toward}, ``The communication infrastructures for implementing distributed methods need to be carefully designed, and the impact of communication delays and failures on the performance of distributed methods need to be investigated''.

Instead of using decomposition or decentralized procedures, which heavily rely on access to either all physical and economic information or fail-safe real-time communication infrastructure, what we propose in this paper is an approximate method that requires neither of these two controversial assumptions at the expense of obtaining a solution slightly different from the optimal one. The proposed approach only needs access to offline historical information on prices and power injection at the substations connecting the transmission and distribution networks. Using statistical tools, we learn the price response of active distribution networks whose operating decisions aim at minimizing costs while complying with local physical constraints, such as voltage limits. We also utilize easily accessible information (e.g., capacity factors of wind and solar local resources) to make the curve adaptive to changes in external conditions that affect power system operations. Finally, the obtained response is approximated by a non-increasing step-wise function that can be conveniently interpreted as a market bid for the participation of the distribution networks in wholesale electricity markets. In summary, the contributions of our paper are twofold:
\begin{enumerate}
    \item We propose a TSO-DSO coordination scheme that uses historical data at substations to learn the price response of active distribution networks using a decreasing step-wise function. If compared with existing methodologies, ours is simple and easy to implement within current market procedures, and cheap in terms of computational resources, information exchange and communication infrastructure. 
    \item We measure the performance of the proposed approach in terms of the power imbalances and social welfare loss caused by the approximation of the distribution networks' behavior using a realistic case study.
\end{enumerate}

We compare our approach against a fully centralized operational model, referred to as \emph{benchmark}, that guarantees the optimal coordination between the TSO and the DSOs. Since this benchmark produces the same solution  obtained with exact decomposition and decentralized methods \cite{kargarian2018toward,yuan_2017,bragin_2017,nawaz_2020,mohammadi2019diagonal}, these have not been considered in our study. On the contrary, our model is evaluated against  two other approximations. In the first one, called the \emph{single-bus approach}, all physical constraints of the distribution networks are disregarded as if all small consumers and distributed generating resources were directly connected to the main substation. In the second one,  called the \emph{price-agnostic approach}, the response of the distribution networks is assumed to be independent of prices.


We note that, within the context of reactive power optimization for the minimization of network losses, the authors in \cite{Ding2018} also approximate the apparent power exchange between the TSO and the DSOs by a polynomial function of the voltage level at the main substation. However, beyond the evident facts that their purpose is different and the fitting procedure we need to use is more intricate (to comply with market rules), they also omit the \emph{dynamic} nature of distribution network response to local marginal prices (LMPs). Similarly, authors of \cite{li2018response} also propose a methodology to obtain the relation between the distribution network response and the voltage level at the substation. However, their approach is based on perfect knowledge of all distribution network parameters.

The rest of this paper is organized as follows. Section~\ref{sec:MF} introduces optimization models for transmission and distribution network operations, which are then used to construct  different  DSO-TSO coordination approaches  in Section~\ref{sec:methodology}. The metrics we use for comparing these approaches are described in Section~\ref{sec:CP}, while the case study is presented in Section~\ref{sec:SR}. Finally, conclusions are duly reported in Section~\ref{sec:conc}.

\section{Modeling Framework}\label{sec:MF}

\begin{figure}
\centering
\begin{tikzpicture}[scale=0.7, every node/.style={scale=0.7}]
\node at (-0.4,-0.7) {\includegraphics[scale=0.03]{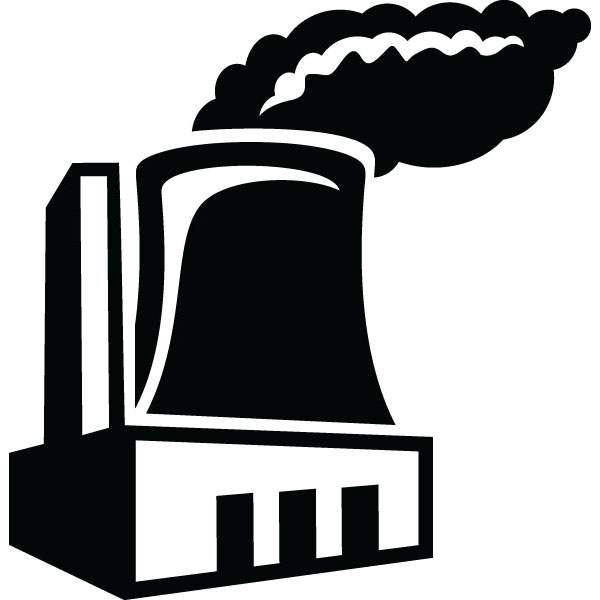}};
\node at (-1.3,0.8) {\includegraphics[scale=0.03]{pictures/powerplant.jpg}};
\node at (1.2,-0.8) {\includegraphics[scale=0.07]{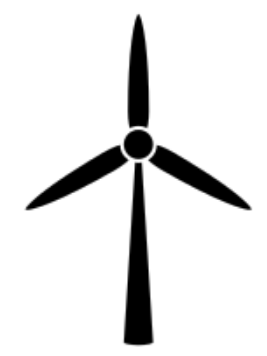}};
\node at (1.5,0.7) {\includegraphics[scale=0.07]{pictures/wind2.png}};
\node at (-4,0.85) {\includegraphics[scale=0.05]{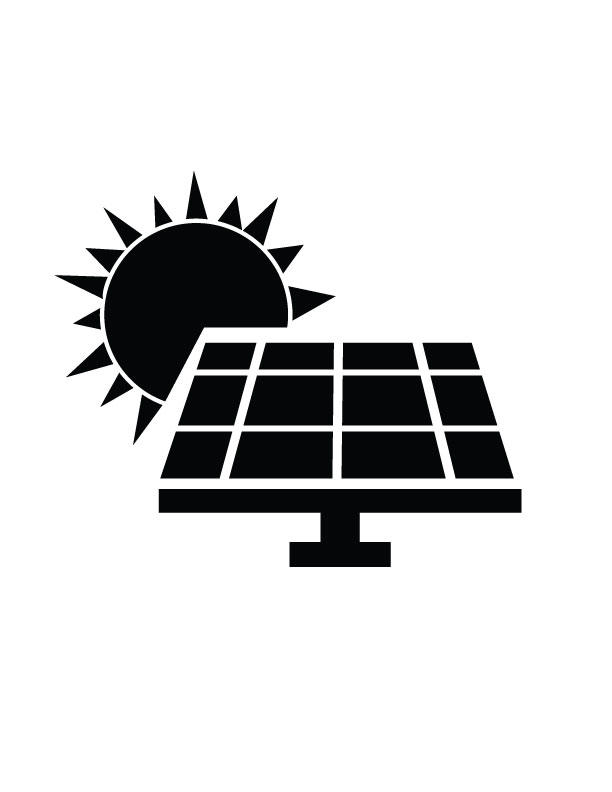}};
\node at (-4.5,-0.6) {\includegraphics[scale=0.08]{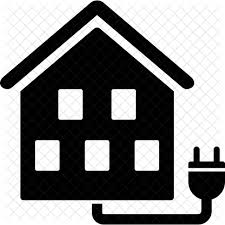}};
\node at (-4+1.5,0.85-3) {\includegraphics[scale=0.05]{pictures/solar.jpg}};
\node at (-4.5+1.5,-0.6-3) {\includegraphics[scale=0.08]{pictures/home.jpeg}};
\node at (4,0.85) {\includegraphics[scale=0.05]{pictures/solar.jpg}};
\node at (4.5,-0.6) {\includegraphics[scale=0.08]{pictures/home.jpeg}};
\node at (4-1.5,0.85-3) {\includegraphics[scale=0.05]{pictures/solar.jpg}};
\node at (4.5-1.5,-0.6-3) {\includegraphics[scale=0.08]{pictures/home.jpeg}};
\node[circle,draw=black, inner
sep=0pt,minimum size=4cm,label=above:{Transmission network}] (c1) at (0,0) {};
\node[circle,draw=black, inner
sep=0pt,minimum size=2.5cm, label=below:{Distribution networks}] (c2) at (4,0) {};
\node[circle,draw=black, inner
sep=0pt,minimum size=2.5cm, label=below:{Distribution networks}] (c3) at (-4,0) {};
\node[circle,draw=black, inner
sep=0pt,minimum size=2.5cm] (c4) at (-2.5,-3) {};
\node[circle,draw=black, inner
sep=0pt,minimum size=2.5cm] (c5) at (2.5,-3) {};
\node[circle,fill,inner sep=1pt] (n1) at (0.6,1.5) {};
\node[circle,fill,inner sep=1pt] (n2) at (0.3,-1.5){};
\node[circle,fill,inner sep=1pt] (n3) at (-0.5,1.2){};
\node[circle,fill,inner sep=1pt] (n4) at (0,-0.2){};
\node[circle,fill,inner sep=1pt] (n5) at (-0.6,-1.5){};
\node[circle,fill,inner sep=1pt] (n6) at (-1,0.4){};
\node[circle,fill,inner sep=1pt] (n7) at (0.6,0.5){};
\node[circle,fill,inner sep=1pt] (n8) at (1.4,0){};
\node[circle,fill,inner sep=1pt] (n9) at (-1.5,0){};
\draw [-] (n1) -- (n3);
\draw [-] (n1) -- (n7);
\draw [-] (n1) -- (n8);
\draw [-] (n3) -- (n7);
\draw [-] (n7) -- (n8);
\draw [-] (n4) -- (n2);
\draw [-] (n6) -- (n4);
\draw [-] (n4) -- (n9);
\draw [-] (n3) -- (n6);
\draw [-] (n4) -- (n7);
\draw [-] (n4) -- (n8);
\draw [-] (n2) -- (n8);
\draw [-] (n2) -- (n5);
\draw [-] (n5) -- (n9);
\draw [-] (n6) -- (n9);
\draw [-] (n1) -- (n3);
\draw [-] (n1) -- (n3);
\node[circle,fill,inner sep=1pt] (n11) at (-3,0){};
\node[circle,fill,inner sep=1pt] (n12) at (-3.5,0.5){};
\node[circle,fill,inner sep=1pt] (n13) at (-3.8,0.5){};
\node[circle,fill,inner sep=1pt] (n14) at (-4.1,0.5){};
\node[circle,fill,inner sep=1pt] (n15) at (-4.4,0.5){};
\node[circle,fill,inner sep=1pt] (n16) at (-3.5,0){};
\node[circle,fill,inner sep=1pt] (n17) at (-3.8,0){};
\node[circle,fill,inner sep=1pt] (n18) at (-4.1,0){};
\node[circle,fill,inner sep=1pt] (n19) at (-4.4,0){};
\node[circle,fill,inner sep=1pt] (n20) at (-4.7,0){};
\node[circle,fill,inner sep=1pt] (n21) at (-5,0){};
\node[circle,fill,inner sep=1pt] (n22) at (-3.5,-0.5){};
\node[circle,fill,inner sep=1pt] (n23) at (-3.8,-0.5){};
\node[circle,fill,inner sep=1pt] (n24) at (-4.1,-0.5){};
\draw [-] (n11) -- (n12);
\draw [-] (n12) -- (n13);
\draw [-] (n13) -- (n14);
\draw [-] (n14) -- (n15);
\draw [-] (n11) -- (n16);
\draw [-] (n16) -- (n17);
\draw [-] (n17) -- (n18);
\draw [-] (n18) -- (n19);
\draw [-] (n19) -- (n20);
\draw [-] (n20) -- (n21);
\draw [-] (n11) -- (n22);
\draw [-] (n22) -- (n23);
\draw [-] (n22) -- (n24);
\draw [-] (n9) -- (n11);
\node[circle,fill,inner sep=1pt] (n31) at (-3+1.5,0-3){};
\node[circle,fill,inner sep=1pt] (n32) at (-3.5+1.5,0.5-3){};
\node[circle,fill,inner sep=1pt] (n33) at (-3.8+1.5,0.5-3){};
\node[circle,fill,inner sep=1pt] (n34) at (-4.1+1.5,0.5-3){};
\node[circle,fill,inner sep=1pt] (n35) at (-4.4+1.5,0.5-3){};
\node[circle,fill,inner sep=1pt] (n36) at (-3.5+1.5,0-3){};
\node[circle,fill,inner sep=1pt] (n37) at (-3.8+1.5,0-3){};
\node[circle,fill,inner sep=1pt] (n38) at (-4.1+1.5,0-3){};
\node[circle,fill,inner sep=1pt] (n39) at (-4.4+1.5,0-3){};
\node[circle,fill,inner sep=1pt] (n40) at (-4.7+1.5,0-3){};
\node[circle,fill,inner sep=1pt] (n41) at (-5+1.5,0-3){};
\node[circle,fill,inner sep=1pt] (n42) at (-3.5+1.5,-0.5-3){};
\node[circle,fill,inner sep=1pt] (n43) at (-3.8+1.5,-0.5-3){};
\node[circle,fill,inner sep=1pt] (n44) at (-4.1+1.5,-0.5-3){};
\draw [-] (n31) -- (n32);
\draw [-] (n32) -- (n33);
\draw [-] (n33) -- (n34);
\draw [-] (n34) -- (n35);
\draw [-] (n31) -- (n36);
\draw [-] (n36) -- (n37);
\draw [-] (n37) -- (n38);
\draw [-] (n38) -- (n39);
\draw [-] (n39) -- (n40);
\draw [-] (n40) -- (n41);
\draw [-] (n31) -- (n42);
\draw [-] (n42) -- (n43);
\draw [-] (n42) -- (n44);
\draw [-] (n5) -- (n31);
\node[circle,fill,inner sep=1pt] (n41) at (3,0){};
\node[circle,fill,inner sep=1pt] (n42) at (3.5,0.5){};
\node[circle,fill,inner sep=1pt] (n43) at (3.8,0.5){};
\node[circle,fill,inner sep=1pt] (n44) at (4.1,0.5){};
\node[circle,fill,inner sep=1pt] (n45) at (4.4,0.5){};
\node[circle,fill,inner sep=1pt] (n46) at (3.5,0){};
\node[circle,fill,inner sep=1pt] (n47) at (3.8,0){};
\node[circle,fill,inner sep=1pt] (n48) at (4.1,0){};
\node[circle,fill,inner sep=1pt] (n49) at (4.4,0){};
\node[circle,fill,inner sep=1pt] (n50) at (4.7,0){};
\node[circle,fill,inner sep=1pt] (n51) at (5,0){};
\node[circle,fill,inner sep=1pt] (n52) at (3.5,-0.5){};
\node[circle,fill,inner sep=1pt] (n53) at (3.8,-0.5){};
\node[circle,fill,inner sep=1pt] (n54) at (4.1,-0.5){};
\draw [-] (n41) -- (n42);
\draw [-] (n42) -- (n43);
\draw [-] (n43) -- (n44);
\draw [-] (n44) -- (n45);
\draw [-] (n41) -- (n46);
\draw [-] (n46) -- (n47);
\draw [-] (n47) -- (n48);
\draw [-] (n48) -- (n49);
\draw [-] (n49) -- (n50);
\draw [-] (n50) -- (n51);
\draw [-] (n41) -- (n52);
\draw [-] (n52) -- (n53);
\draw [-] (n52) -- (n54);
\draw [-] (n8) -- (n41);
\node[circle,fill,inner sep=1pt] (n61) at (3-1.5,0-3){};
\node[circle,fill,inner sep=1pt] (n62) at (3.5-1.5,0.5-3){};
\node[circle,fill,inner sep=1pt] (n63) at (3.8-1.5,0.5-3){};
\node[circle,fill,inner sep=1pt] (n64) at (4.1-1.5,0.5-3){};
\node[circle,fill,inner sep=1pt] (n65) at (4.4-1.5,0.5-3){};
\node[circle,fill,inner sep=1pt] (n66) at (3.5-1.5,0-3){};
\node[circle,fill,inner sep=1pt] (n67) at (3.8-1.5,0-3){};
\node[circle,fill,inner sep=1pt] (n68) at (4.1-1.5,0-3){};
\node[circle,fill,inner sep=1pt] (n69) at (4.4-1.5,0-3){};
\node[circle,fill,inner sep=1pt] (n70) at (4.7-1.5,0-3){};
\node[circle,fill,inner sep=1pt] (n71) at (5-1.5,0-3){};
\node[circle,fill,inner sep=1pt] (n72) at (3.5-1.5,-0.5-3){};
\node[circle,fill,inner sep=1pt] (n73) at (3.8-1.5,-0.5-3){};
\node[circle,fill,inner sep=1pt] (n74) at (4.1-1.5,-0.5-3){};
\draw [-] (n61) -- (n62);
\draw [-] (n62) -- (n63);
\draw [-] (n63) -- (n64);
\draw [-] (n64) -- (n65);
\draw [-] (n61) -- (n66);
\draw [-] (n66) -- (n67);
\draw [-] (n67) -- (n68);
\draw [-] (n68) -- (n69);
\draw [-] (n69) -- (n70);
\draw [-] (n70) -- (n71);
\draw [-] (n61) -- (n72);
\draw [-] (n72) -- (n73);
\draw [-] (n72) -- (n74);
\draw [-] (n2) -- (n61);
\end{tikzpicture}
\caption{Transmission and distribution coordination scheme}
\label{fig:tso_dso}
\end{figure}

We consider a power system with a high-voltage, meshed transmission network  connected to generating units, large consumers  and several medium-voltage distribution networks. As illustrated in Fig. \ref{fig:tso_dso}, each distribution system is connected to the transmission network through one main substation,  has a radial topology and hosts small-scale electricity consumers and producers. 

The active power output of generating unit $i$ at time period $t$ is denoted by $\pg$, with  minimum/maximum limits $\underline{p}^G_i/\overline{p}^G_i$. Generating units are assumed to have a convex cost function of the form $c_i(\pg)=\frac{1}{2}a_i(\pg)^2 + b_i(\pg)$, with $a_i,b_i \geq 0$, and a dimensionless capacity factor $\rho_{it}$, with $0\leq\rho_{it}\leq1$. For thermal units $\rho_{it}=1, \forall t$, while for renewable generating units the capacity factor depends on weather conditions and the production cost is zero ($a_i=b_i=0$).

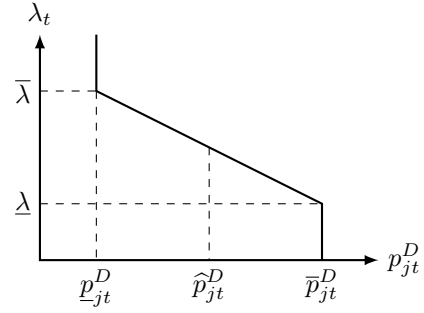
\begin{figure}
    \centering
    \begin{tikzpicture}[scale=1.5]
    \draw [latex-latex,thick] (0,2) node (yaxis) [above] {$\lambda_t$}
        |- (3,0) node (xaxis) [right] {$p^D_{jt}$};
    \draw[thick] (0.5,2) -- (0.5,1.5) -- (2.5,0.5) -- (2.5,0);
    \draw[dashed] (0,1.5) --(0.5,1.5);
    \draw[dashed] (0,0.5)--(2.5,0.5);
    \draw[dashed] (1.5,1)--(1.5,0);
    \draw[dashed] (0.5,1.5)--(0.5,0);
    \node[below] at (1.5,0) {$\widehat{p}^D_{jt}$};
    \node[below] at (0.5,0) {$\underline{p}^D_{jt}$};
    \node[below] at (2.5,0) {$\overline{p}^D_{jt}$};
    \node[left] at (0,1.5) {$\overline{\lambda}$};
    \node[left] at (0,0.5) {$\underline{\lambda}$};
\end{tikzpicture}
    \caption{Flexible electricity demand modeling}
    \label{fig:flex}
\end{figure}

Electricity consumption is modeled as a capped linear function of the LMP $\lambda_t$, as shown in Fig. \ref{fig:flex}, where $\widehat{p}^D_{jt}$ denotes the baseline demand of consumer $j$ at time $t$ and $\overline{p}^D_{jt}/ \underline{p}^D_{jt}$ are the maximum/minimum load levels given by $\overline{p}^D_{jt}=\widehat{p}^D_{jt}(1+\delta_j)$ and $\underline{p}^D_{jt}=\widehat{p}^D_{jt}(1-\delta_j)$, with $\delta_j \geq 0$ \cite{Mieth2020}. Consequently, the maximum/minimum load levels vary over time according to the evolution of the baseline demand. Under this modeling approach, a price-insensitive demand is modeled with $\delta_j=0$, while $\delta_j=0.5$ implies that the consumer is willing to increase or decrease their baseline demand up to 50\% depending on the price.  Finally, $\overline{\lambda}$ and $\underline{\lambda}$ stand for the LMP values that unlock the minimum and maximum demand from consumers, respectively. The  demand function in Fig. \ref{fig:flex} goes through points $(\underline{p}^D_{jt},\overline{\lambda})$ and $(\overline{p}^D_{jt},\underline{\lambda})$ and, therefore, its expression can be determined as follows:
\begin{align}
& \frac{\lambda_t - \underline{\lambda}}{\overline{\lambda}-\underline{\lambda}} = \frac{p^D_{jt}-\overline{p}^D_{jt}}{\underline{p}^D_{jt}-\overline{p}^D_{jt}} = \frac{p^D_{jt}-\widehat{p}^D_{jt}(1+\delta_j)}{\widehat{p}^D_{jt}(1-\delta_j)-\widehat{p}^D_{jt}(1+\delta_j)} \implies \nonumber \\
& p^D_{jt} = \widehat{p}^D_{jt}\left( 1 + \delta_j \frac{\overline{\lambda}+\underline{\lambda}}{\overline{\lambda}-\underline{\lambda}} \right) - \frac{2\widehat{p}^D_{jt}\delta_j}{\overline{\lambda}-\underline{\lambda}} \lambda_t 
\end{align}

Hence, the active demand level $\pd$ for a given electricity price $\lambda_t$ takes the following form:
\begin{equation}
\pd = \left\{ \begin{array}{lcl}
\overline{p}^D_{jt} &  \text{if}  & \lambda_t  \leq \underline{\lambda} \\
\alpha_{jt} - \beta_{jt} \lambda_t &  \text{if} & \underline{\lambda} < \lambda_t < \overline{\lambda}
\\ 
\underline{p}^D_{jt} &   \text{if}  &  \overline{\lambda} \leq \lambda_t,
\end{array}
\right.   \label{eq:inverse_function} 
\end{equation}
where $\alpha_{jt}=\widehat{p}^D_{jt}\left(1+\delta_j\frac{\overline{\lambda}+\underline{\lambda}}{\overline{\lambda}-\underline{\lambda}} \right)$ and $\beta_{jt}=\frac{2\widehat{p}^D_{jt}\delta_j}{\overline{\lambda}-\underline{\lambda}}$. The reactive power demand is given by $\qd=\gamma_j\pd$, where $\gamma_j$ is the power factor of consumer $j$, which is assumed to be independent of time for simplicity. Finally, we obtain the utility of each consumer by integrating the inverse demand function with respect to the demand quantity, that is,
\begin{align}
& u_{jt}(p^D_{jt}) = \int_{\underline{p}^D_{jt}}^{p^D_{jt}} \lambda(p) dp = \int_{\underline{p}^D_{jt}}^{p^D_{jt}} \left( \frac{\alpha_{jt}}{\beta_{jt}} - \frac{p}{\beta_{jt}} \right) dp = \nonumber \\
& \qquad = \frac{\alpha_{jt}}{\beta_{jt}} \left(p^D_{jt}-\underline{p}^D_{jt}\right) - \frac{(p^D_{jt})^2-(\underline{p}^D_{jt})^2}{2\beta_{jt}} \label{eq:utility}
\end{align}
%

The transmission network is modeled using a DC power flow approximation \cite{nawaz2020stochastically} and, therefore, each line $l$ going from node $o_l$ to node $e_l$ is characterized by its reactance $x_l$ and maximum capacity $\overline{s}^F_l$. The power flow is denoted by $\pf$.

Then, suppose that the active consumption of the $k$-th distribution network $p^N_{kt}$ can be expressed as $p^N_{kt}=h_{kt}(\lambda_{kt})$, where $\lambda_{kt}$ is the price at the corresponding substation. Under this assumption, transmission system operations at time period $t$ are modeled by the following optimization problem:
\begin{subequations}
\begin{align}
    & \max_{\Phi^T_t} \sum_{j\in J^T} u_{jt}(\pd) + \sum_{k\in K^T} \int_{0}^{p^N_{kt}} \hspace{-3mm} h_{kt}^{-1}(s)ds - \sum_{i\in I^T} c_i(\pg) 
    \label{eq:transmission_of}\\
    & \text{s.t.} \nonumber \\ 
    & \sum_{i\in I_n} \pg - \sum_{j\in J_n} \pd - \sum_{k\in K_n} p^N_{kt} = \nonumber \\
    & \qquad = \sum_{l:e_l=n} \pf - \sum_{l:o_l=n} \pf, \; \forall n\in N^T \label{eq:transmission_bal} \\
    & \frac{\pf}{s^B} = \frac{1}{x_l}(\theta_{o_lt}-\theta_{e_lt}), \;\forall l\in L^T \label{eq:transmission_flow} \\
    & \underline{p}^G_i \leq \pg \leq \rho_{it}\overline{p}^G_i, \;\forall i\in I^T \label{eq:transmission_maxgen}\\
    & \underline{p}^D_{jt} \leq \pd \leq \overline{p}^D_{jt}, \;\forall j\in J^T \label{eq:transmission_maxdem}\\
    & - \overline{s}^F_l \leq \pf \leq \overline{s}^F_l, \;\forall l\in L^T \label{eq:transmission_maxflow}
\end{align} \label{eq:transmission}
\end{subequations}
\noindent where $\theta_{nt}$ is the voltage angle at node $n$ and time period $t$, $\Phi^T_t=(\pg,\pd,p^N_{kt},\pf,\theta_{nt})$ are decisions variables,  $N^T,L^T,I^T,J^T,K^T$ are sets of  nodes, lines, generators, consumers and distribution networks connected to the transmission network, and $I_n,J_n,K_n$ are sets  of generating units, consumers and distribution networks connected to node $n$. Objective function \eqref{eq:transmission_of} maximizes the total social welfare and includes the utility of all flexible consumers connected to the transmission network (first term), the utility of all distribution networks (second term), and the generation cost of all units connected to the transmission network (third term). Note that  $h_{kt}^{-1}(\cdot)$ represents the inverse demand function and its integral correspond to the total utility of each distribution network. The nodal power balance equation is imposed by \eqref{eq:transmission_bal}, while the power flow through each transmission line is computed in \eqref{eq:transmission_flow}. Finally, constraints \eqref{eq:transmission_maxgen},  \eqref{eq:transmission_maxdem} and \eqref{eq:transmission_maxflow}  enforce the generation, consumption and transmission capacity limits. 

Traditionally, distribution networks only hosted inflexible consumption and, therefore, $p^N_{kt}$ was considered independent of the electricity price. In this case, the second term of \eqref{eq:transmission_of} vanishes, and  variable $p^N_{kt}$ is replaced by the forecast power intake of each distribution network. Thus, problem \eqref{eq:transmission} can be transformed into a quadratic optimization problem that can be solved to global optimality using off-the-shelf solvers, \cite[Appendix B]{doi:10.1137/1.9781611974164.fm}. However, this paradigm has changed in recent years and current distribution networks include a growing amount of flexible small-scale consumers and distributed generation resources that are capable of  adjusting their consumption/generation in response to the electricity price to maximize their utility/payoff \cite{Papavasiliou2018}. Indeed, if $\lambda_{kt}$ is the electricity price at the main substation of distribution network $k$, the power intake of that distribution network $p^N_{kt}$ can be determined by solving the following optimization problem:
\begin{subequations}
\begin{align}
    & \max_{\Phi^D_{kt}} \quad \sum_{j\in J_k^D} u_{jt}(\pd) -   \sum_{i\in I_k^D} c_i(\pg) - \lambda_{kt}p^N_{kt}   \label{eq:distribution_of}\\
    & \text{s.t.} \nonumber \\
    & p^N_{kt} + \sum_{i\in I_n} \pg - \sum_{j\in J_n} \pd = \nonumber \\
    & \qquad = \sum_{l:e_l=n} \pf - \sum_{l:o_l=n} \pf, \; n=n^0_k \label{eq:distribution_bal0}\\
    & \sum_{i\in I_n} \pg - \sum_{j\in J_n} \pd = \nonumber \\
    & \qquad = \sum_{l:e_l=n} \pf - \sum_{l:o_l=n} \pf, \; \forall n \in N_k^D, n\neq n^0_k \label{eq:distribution_bal}\\
    & \sum_{i\in I_n} \qg - \sum_{j\in J_n} \qd = \nonumber \\
    & \qquad = \sum_{l:e_l=n} \qf - \sum_{l:o_l=n} \qf, \; \forall n \in N_k^D\label{eq:distribution_balq}\\
    & \qd = \gamma_j \pd, \; \forall  j \in J_k^D \label{eq:distribution_factor}\\
    & v_{nt} = v_{a_nt} - \frac{2}{s^B} \sum_{l:e_l=n} r_l\pf + x_l\qf, \; \forall n \in N_k^D \label{eq:distribution_voltage}\\
    & \underline{p}^G_i \leq \pg \leq \rho_{it}\overline{p}^G_i, \; \forall i \in I_k^D \label{eq:distribution_maxgen}\\
    & \underline{q}^G_i \leq \qg \leq \overline{q}^G_i, \; \forall i \in I_k^D \label{eq:distribution_maxgenq}\\
    & (\pg)^2 + (\qg)^2 \leq (\overline{s}^G_i)^2, \; \forall i \in I_k^D \label{eq:distribution_converter}\\
    & \underline{p}^D_{jt} \leq \pd \leq \overline{p}^D_{jt}, \; \forall  j \in J_k^D \label{eq:distribution_maxdem}\\
    & (\pf)^2 + (\qf)^2 \leq (\overline{s}^F_l)^2, \; \forall l \in L_k^D \label{eq:distribution_maxflow}\\
    & \underline{v}_{nt} \leq v_{nt} \leq \overline{v}_{nt}, \; \forall n \in N_k^D \label{eq:distribution_maxvol}
\end{align} \label{eq:distribution}
\end{subequations}
where the decisions variables are $\Phi^D_{kt} = (p^N_{kt}, \pg, \qg, \pd, \qd, \pf, \qf, v_{nt})$. In particular,  $\qg,\qd,\qf$ are the reactive power generation, consumption and flow, in that order, and $v_{nt}$ is the squared voltage magnitude. Since we assume a radial distribution network, we use the LinDistFlow AC power flow approximation, where $a_n$ represents the ancestor of node $n$ and $r_l$ is the resistance of line $l$ \cite{Mieth2018}. The rate power of the inverters for distributed generators is denoted as $\overline{s}^G_i$ \cite{Hassan2018a}, and the squared voltage magnitude limits are  $\underline{v}_{nt},\overline{v}_{nt}$. Finally, 
$N_k^D,L_k^D,I_k^D,J_k^D$ are the set of nodes, lines, generators and consumers of distribution network $k$, and $n^0_k$ corresponds to the node of the distribution network connected to the substation. 

Objective function \eqref{eq:distribution_of} maximizes the social welfare of distribution network $k$ and includes the utility of flexible consumers (first term), the cost of distributed generation (second term) and the cost of power exchanges with the transmission network (third term). Nodal active and reactive power equations are formulated in~\eqref{eq:distribution_bal0}, \eqref{eq:distribution_bal} and \eqref{eq:distribution_balq}. Constraint \eqref{eq:distribution_factor} relates active and reactive demand through a given power factor, while the dependence of voltage magnitudes in a radial network is accounted for in \eqref{eq:distribution_voltage} using the LinDistFlow approximation. Limits on active and reactive generating power outputs are enforced in \eqref{eq:distribution_maxgen}, \eqref{eq:distribution_maxgenq} and \eqref{eq:distribution_converter}. Similarly, equations \eqref{eq:distribution_maxdem}, \eqref{eq:distribution_maxflow} and \eqref{eq:distribution_maxvol} determine the feasible values of demand quantities, power flows and squared voltage magnitudes. As a result,  \eqref{eq:distribution} is a convex optimization problem that can be solved using off-the-shelf solvers. 

Drawing a closed-form expression $h_{kt}(\lambda_{kt})$ from~\eqref{eq:distribution} that exactly characterizes the optimal value of $p^N_{kt}$ as a function of the electricity price $\lambda_{kt}$ seems like a lost cause. Furthermore, even if such an expression were possible, using it in \eqref{eq:transmission} would lead to a troublesome  non-convex optimization problem, with the likely loss of global optimality guarantees. In the next section, we discuss different strategies to construct an approximation  $\hat{h}_{kt}(\lambda_{kt})$ that can be easily incorporated into \eqref{eq:transmission} to determine the optimal operation of the transmission network. In particular, we focus on strategies that leverage available contextual information to construct function $\hat{h}_{kt}(\lambda_{kt})$.

\section{Methodology} \label{sec:methodology}

In this section we present four different approaches to accommodate the behavior of active distribution networks in transmission network operations. The aim of these methods is to determine the electricity prices at substations that foster the most efficient use of the flexible resources available in the distribution networks. Once electricity prices are published, the DSOs operate the distributed energy resources to maximize their social welfare while satisfying the physical limits of the distribution network, such as voltage limits and reactive power capacities.

\subsection{Benchmark approach (BN)}

This approach includes a full representation of both the transmission system and the distribution networks, by jointly solving optimization problems \eqref{eq:transmission} and \eqref{eq:distribution} as follows:
\begin{subequations}
\begin{align}
& \max_{\Phi^T_t,\Phi^D_{kt}} \sum_{j\in J^T \cup \{J^D_k\}} u_{jt}(\pd) - \sum_{i\in I^T \cup \{I^D_k\}} c_i(\pg)  \label{eq:benchmark_of}\\  
    & \text{s.t.} \qquad \eqref{eq:transmission_bal}-\eqref{eq:transmission_maxflow}, \eqref{eq:distribution_bal}-\eqref{eq:distribution_maxvol} 
\end{align} \label{eq:benchmark}
\end{subequations}
Model \eqref{eq:benchmark} enables the optimal operation of the transmission network since it takes into account the most accurate representation of all distribution networks connected to it \cite{LeCadre2019}. However, this approach has the following drawbacks: 
\begin{itemize}
    \item[-] It requires having access to distribution network parameters, such as its topological configuration and $r_l,x_l$, which is impractical, as private  or sovereign entities operating distribution networks prefer to keep this information confidential \cite{mohammadi2019diagonal,Ding2018,Yu2018}.
    \item[-] Operating the power system through \eqref{eq:benchmark} would require a deep transformation of current market mechanisms to allow small generators/consumers to directly submit their electricity offers/bids to a centralized market operator. 
    \item[-] Even if all distribution network parameters were known and small generators/consumers were allowed to directly participate in the electricity market, solving model \eqref{eq:benchmark} is computationally expensive for realistically sized systems with hundreds of distribution networks connected to the transmission network \cite{Li2018}.
\end{itemize}

In this paper, we use the solution of this approach as a benchmark to evaluate the performance of the other methods described in this section. Other decomposition or decentralized methods in the technical literature are able to achieve global optimality and, therefore, their solution coincide with that of BN. For this reason, we focus on comparing the proposed approach with other approximate methodologies that also lead to suboptimal solutions.

\subsection{Single-bus approach (SB)}

This approach is a relaxation of BN in  \eqref{eq:benchmark}, where  physical limits on distribution power flows and voltages are disregarded. Therefore, operational model SB can be equivalently interpreted as if all small consumers and distributed energy resources were directly connected to the transmission network, i.e. all  distribution systems are modeled as single-bus grids. Therefore, the dispatch decisions for the transmission network are computed by solving the following problem:
\begin{subequations}
\begin{align}
& \max_{\Phi^T_t,\Phi^D_{kt}} \sum_{j\in J^T \cup \{J^D_k\}} u_{jt}(\pd) - \sum_{i\in I^T \cup \{I^D_k\}} c_i(\pg) \label{eq:dso2tso_of}\\ 
& \text{s.t.} \nonumber \\
& \sum_{i\in \hat{G}_n} \pg - \sum_{j\in \hat{D}_n} \pd = \sum_{l:e_l=n} \pf - \sum_{l:o_l=n} \pf, \;\forall n\in N^T \label{eq:dso2tso_bal} \\
& \frac{\pf}{s^B} = \frac{1}{x_l}(\theta_{o_lt}-\theta_{e_lt}), \;\forall l\in L^T \label{eq:dso2tso_flow} \\
& \underline{p}^G_i \leq \pg \leq \rho_{it}\overline{p}^G_i, \;\forall i\in I^T \cup \{I^D_k\} \label{eq:dso2tso_maxgen}\\
& \underline{p}^D_{jt} \leq \pd \leq \overline{p}^D_{jt}, \;\forall j\in J^T \cup \{J^D_k\} \label{eq:dso2tso_maxdem}\\
& - \overline{s}^F_l \leq \pf \leq \overline{s}^F_l, \;\forall l\in L^T \label{eq:dso2tso_maxflow}
\end{align} \label{eq:dso2tso}
\end{subequations}
where $\hat{G}_n$ and $\hat{D}_n$ denote, respectively, the set of generators and consumers either directly connected to node $n$ or hosted by a distribution network connected to it. Problem \eqref{eq:dso2tso} is less computationally demanding than the BN approach  in \eqref{eq:benchmark} and does not require knowledge of distribution network parameters. However, this approach also relies on a market mechanism that allows small generators and consumers to submit their offers and bids directly to the wholesale market \cite{Chen2018}. Besides, if the operation of some of the distribution networks is constrained by the physical limitations of power flows and/or voltage levels, then the solution provided by this approach may substantially differ from the actual   conditions in the distribution networks. 

\subsection{Contextual price-agnostic approach (PAG)}

This approach is based on the premise that the penetration rates of small-scale flexible consumers and distributed generation resources is not significant and, therefore, the response of  distribution networks is independent of LMPs at their substations. On the other hand, this response can still depend on other contextual information that affect the behavior of distribution networks such as the aggregated load level of their flexible consumers and the wind and solar capacity factors in the corresponding geographical area. 

Consider a set of historical data $\{\chi_{kt},p^N_{kt}\}_{t \in \mathcal{T}}$, where $\chi_{kt}$ represents a vector containing the contextual information to explain the consumption level of distribution network $k$. Vector $\chi_{kt}$ can include weather conditions, e.g. ambient temperature, wind speed, solar irradiation, or categorical variables, e.g. an hour of the day or a day of the week. The PAG approach aims to learn the relation between $p^N_{kt}$ and $\chi_{kt}$ for each distribution network $k$, i.e.,
\begin{equation}
p^N_{kt} = f_k(\chi_{kt}) \label{eq:contex1}
\end{equation}
The function $f_k$ that best approximates the behavior of distribution network $k$ with contextual information can be found using a wide variety of supervised learning techniques \cite{Hastie2009}. In particular, if $f_k$ must belong to a certain family of functions, such as the family of linear functions, its parameters can be computed using the well-known least squares criterion. In order to capture non-linear relations, $f_k$ can also represent a neural network to be trained using available data. Alternatively, if we do not make strong assumptions about the form of the mapping function, the relation between $\chi_{kt}$ and $p^N_{kt}$ can be modeled using non-parametric supervised learning techniques. Within this group, we opt in this work for the k-nearest neighbors regression algorithm ($K$-NN) because of its simplicity, interpretability and scalability. Following this methodology, the estimation of the power import of distribution network $k$ for time period $t$ (denoted as $\widehat{p}^N_{kt}$) is computed as:
\begin{equation}
\widehat{p}^N_{kt} = \frac{1}{K}\sum_{t'\in \mathcal{T}^C(t)}p^N_{kt'}, \label{eq:knn} 
\end{equation}
where $\mathcal{T}^C(t)$ is the subset of the $K$ time periods whose contexts are the closest to $\chi_{kt}$ according to a given distance, and $t'$ is an auxiliary time period index. If contextual information only includes continuous variables (electricity demand, renewable power generation, etc.), the dissimilarity between two time periods can be measured using the Euclidean distance, i.e., $dist(t_1,t_2)=||\chi_{kt_1}-\chi_{kt_2}||_2$. If contextual information also includes binary variables (equipment status, maintenance schedules, etc.), the dissimilarity can be measured using the Hamming distance, for example.

Once the forecast intake for each distribution network is obtained depending on its corresponding context, we model all distribution networks as fix loads and determine the operation of the transmission network by solving the following optimization problem:
\begin{subequations}
\begin{align}
    & \max_{\Phi^T} \sum_{j\in J^T} u_{jt}(\pd) - \sum_{i\in I^T} c_i(\pg) \label{eq:average_of}\\
    & \text{s.t.} \nonumber \\
    & p^N_{kt} = \widehat{p}^N_{kt}, \;\forall k \in K \\
    & \eqref{eq:transmission_bal}- \eqref{eq:transmission_maxflow}
\end{align} \label{eq:average}
\end{subequations}
Problem \eqref{eq:average} is also more computationally tractable than   \eqref{eq:benchmark} and does not rely on  knowledge of distribution network parameters. As the one we propose, this approach only requires access to historical power flow measurements at the substation and contextual information that can enhance explainability and interpretability of the distribution network responses to external factors of interest. Fortunately, independent system operators such as ISONE and NYISO make this information publicly available. Another advantage of this approach is that, unlike the SB approach, it can be seamlessly implemented in existing market-clearing procedures since the response of distribution networks is simply replaced with the fixed power injections provided by \eqref{eq:knn}. Actually, this is the approach that most closely reproduces the traditional way of proceeding. On the other hand, since the impact of substation LMPs on the response of distribution networks is disregarded, the accuracy of this approach  worsens as the flexibility  provided by small consumers and  distributed generators increase.


\subsection{Contextual price-aware approach (PAW)}

The SB and PAG approaches disregard the impact of either physical limits or economic signals on the response of distribution networks with small-scale, flexible consumers and distributed generation resources. To overcome this drawback, we propose to approximate the response function $h_{kt}(\lambda_{kt})$ by taking into account the effects of both physical and economic conditions on the behavior of active distribution networks.  
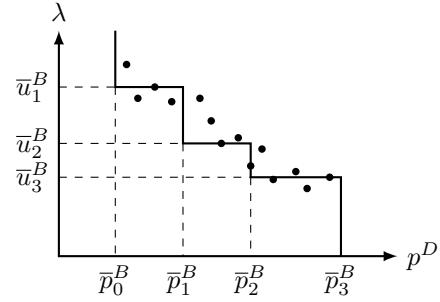
\begin{figure}
    \centering
    \begin{tikzpicture}[scale=1.5]
    \draw [latex-latex,thick] (0,2) node (yaxis) [above] {$\lambda$}
        |- (3,0) node (xaxis) [right] {$p^D$};
    \node at (0.6,1.7) [circle,fill,inner sep=1pt]{};
    \node at (0.7,1.4) [circle,fill,inner sep=1pt]{};
    \node at (0.85,1.5) [circle,fill,inner sep=1pt]{};
    \node at (1.0,1.37) [circle,fill,inner sep=1pt]{};
    \node at (1.25,1.4) [circle,fill,inner sep=1pt]{};
    \node at (1.35,1.2) [circle,fill,inner sep=1pt]{};
    \node at (1.44,1.0) [circle,fill,inner sep=1pt]{};
    \node at (1.59,1.05) [circle,fill,inner sep=1pt]{};
    \node at (1.8,0.95) [circle,fill,inner sep=1pt]{};
    \node at (1.7,0.8) [circle,fill,inner sep=1pt]{};
    \node at (1.9,0.68) [circle,fill,inner sep=1pt]{};
    \node at (2.1,0.75) [circle,fill,inner sep=1pt]{};
    \node at (2.2,0.6) [circle,fill,inner sep=1pt]{};
    \node at (2.4,0.7) [circle,fill,inner sep=1pt]{};
    \draw[thick] (0.5,2) -- (0.5,1.5) -- (1.1,1.5) -- (1.1,1.0) -- (1.7,1.0) -- (1.7,0.7) -- (2.5,0.7) -- (2.5,0);
    \draw[dashed] (0,1.5) --(0.5,1.5);
    \draw[dashed] (0,1.0) --(1.1,1.0);
    \draw[dashed] (0,0.7) --(1.7,0.7);
    \draw[dashed] (0.5,1.5) --(0.5,0);
    \draw[dashed] (1.1,1.0) --(1.1,0);
    \draw[dashed] (1.7,0.7) --(1.7,0);
    \node[below] at (0.5,0) {$\overline{p}^B_0$};
    \node[below] at (1.1,0) {$\overline{p}^B_1$};
    \node[below] at (1.7,0) {$\overline{p}^B_2$};
    \node[below] at (2.5,0) {$\overline{p}^B_3$};
    \node[left] at (0,1.5) {$\overline{u}^B_1$};
    \node[left] at (0,1.0) {$\overline{u}^B_2$};
    \node[left] at (0,0.7) {$\overline{u}^B_3$};
\end{tikzpicture}
    \caption{Step-wise approximation of distribution network response}
    \label{fig:stepwise}
\end{figure}

Similarly to the PAG approach, we assume access to the set of historical data $\{\chi_{kt},\lambda_{kt},p^N_{kt}\}_{t \in \mathcal{T}}$, where $\lambda_{kt}$ denotes the LMP at the substation of distribution network $k$. The proposed PAW approach aims at determining the function that explains the response $p^N_{kt}$ as a function of the contextual information $\chi_{kt}$ and the electricity price at the substation $\lambda_{kt}$, i.e.,
\begin{equation}
p^N_{kt} = g_k(\chi_{kt},\lambda_{kt}) \label{eq:contex2}
\end{equation}
For a fixed context, function \eqref{eq:contex2} provides the relation between the response of a distribution network and the price at its substation. This function can be understood as the bid to be submitted by each distribution network to the wholesale electricity market. However, most current market procedures only accept a finite number of decreasing block bids. For instance, the Spanish market operator establishes that ``For each hourly scheduling period within the same day-ahead scheduling horizon, there can be as many as 25 power blocks for the same production unit, with a different price for each of the said blocks, with the prices increasing for sale bids, or decreasing for purchase bids.'' \cite{OMIE}. In order to comply with these market rules, we propose an efficient learning procedure to determine a decreasing step-wise mapping between the response of a distribution network and the LMP, while taking into account contextual information. Our approach combines unsupervised and supervised  learning techniques as follows:
\begin{itemize}
    \item[-] \textit{Unsupervised learning}. Similarly to PAG, the first step of the proposed approach uses a K-nearest neighbors algorithm to find the subset of time periods $\mathcal{T}^C(t)$ whose contextual information are the closest to $\chi_{kt}$. For the sake of illustration, the points depicted in Fig. \ref{fig:stepwise} represent the pairs of prices and power intakes for times periods in $\mathcal{T}^C(t)$ for a given substation $k$ and context $\chi_{kt}$. 
    \item[-] \textit{Supervised learning}. The second step of the proposed approach consists of finding the step-wise decreasing function that best approximates the price-quantity pairs obtained in the previous step. As illustrated in Fig. \ref{fig:stepwise}, this function can be defined by a set of price breakpoints $\overline{u}^B_{b}$ and the demand level for each block $\overline{p}^B_{b}$. Despite its apparent simplicity, finding the optimal step-wise decreasing function that approximates a set of data points is a complex task that cannot be accomplished by conventional regression techniques. For instance, isotonic regression yields a monotone step-wise function, but a maximum number of blocks cannot be imposed. Conversely, segmented regression provides a step-wise function with a maximum number of blocks, but monotonicity is not ensured. Therefore, the statistical estimation of $\overline{p}^B_{0}$ and  $\overline{u}^B_{b}$ and $\overline{p}^B_{b}$, $\forall b \in \mathcal{B}$, is conducted by means of the curve-fitting algorithm for segmented isotonic regression that has been recently developed in \cite{Bucarey2020} and can be formulated as the following optimization problem:
\begin{subequations} \label{eq:isotonic}
\begin{align}
\hspace{-3mm} \min_{\overline{p}^B_{0},\overline{u}^B_{b},\overline{p}^B_{b}} &  \sum_{t'\in\mathcal{T}^C(t)} {\left(p^N_{kt'}- \hspace{-4mm} \sum_{b \in  \mathcal{B} \cup \{0\}  }{\hspace{-2mm}\overline{p}^B_b \mathbb{I}_{[\overline{u}^B_{b+1}, \overline{u}^B_b)}   (\lambda_{kt'})}\right)^2} \hspace{-3mm} \label{eq:isotonic_of}\\
\textrm{s.t.} & \enskip \overline{p}^B_b \geq \overline{p}^B_{b-1}, \enskip \forall b \in B \label{eq:isotonic_c1}\\
\phantom{\textrm{s.t.}}& \enskip \overline{u}^B_{b+1} \leq \overline{u}^B_{b}, \enskip  \forall b \in B \label{eq:isotonic_c2}
\end{align}
\end{subequations}
\noindent where  $\mathbb{I}_{[\overline{u}^B_{b+1},\overline{u}^B_{b})}$ is the indicator function equal to 1 if $\overline{u}^B_{b+1} \leq \lambda_{kt} < \overline{u}^B_{b}$, and 0 otherwise, and $\overline{u}^B_{0} = \infty$, $\overline{u}^B_{|\mathcal{B}|+1} =-\infty$. Objective function \eqref{eq:isotonic_of} minimizes the sum of squared errors, while constraints \eqref{eq:isotonic_c1}-\eqref{eq:isotonic_c2} ensures the monotonicity of the regression function. Problem \eqref{eq:isotonic} can be reformulated as a mixed-integer quadratic problem to be solved by standard optimization solvers. However, the computational burden of this solution strategy is extremely high. Alternatively, reference \cite{Bucarey2020} proposes a dynamic programming reformulation that guarantees global optimality in polynomial time, which makes this approach computationally attractive.  
\end{itemize}

\begin{figure}
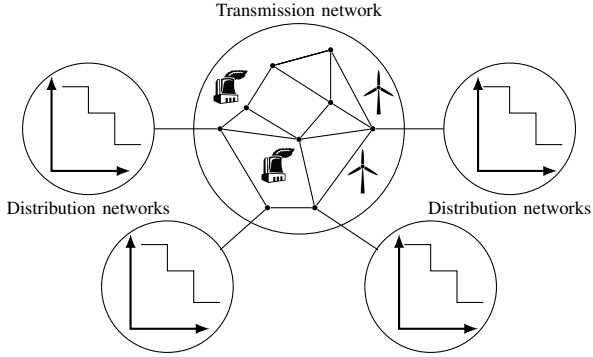

    \centering
    \begin{tikzpicture}[scale=0.7, every node/.style={scale=0.7}]
\node at (-0.4,-0.7) {\includegraphics[scale=0.03]{pictures/powerplant.jpg}};
\node at (-1.3,0.8) {\includegraphics[scale=0.03]{pictures/powerplant.jpg}};
\node at (1.2,-0.8) {\includegraphics[scale=0.07]{pictures/wind2.png}};
\node at (1.5,0.7) {\includegraphics[scale=0.07]{pictures/wind2.png}};
\node[circle,draw=black, inner
sep=0pt,minimum size=4cm,label=above:{Transmission network}] (c1) at (0,0) {};
\node[circle,draw=black, inner
sep=0pt,minimum size=2.5cm, label=below:{Distribution networks}] (c2) at (4,0) {};
\node[circle,draw=black, inner
sep=0pt,minimum size=2.5cm, label=below:{Distribution networks}] (c3) at (-4,0) {};
\node[circle,draw=black, inner
sep=0pt,minimum size=2.5cm] (c4) at (-2.5,-3) {};
\node[circle,draw=black, inner
sep=0pt,minimum size=2.5cm] (c5) at (2.5,-3) {};
\node[circle,fill,inner sep=1pt] (n1) at (0.6,1.5){};
\node[circle,fill,inner sep=1pt] (n2) at (0.3,-1.5){};
\node[circle,fill,inner sep=1pt] (n3) at (-0.5,1.2){};
\node[circle,fill,inner sep=1pt] (n4) at (0,-0.2){};
\node[circle,fill,inner sep=1pt] (n5) at (-0.6,-1.5){};
\node[circle,fill,inner sep=1pt] (n6) at (-1,0.4){};
\node[circle,fill,inner sep=1pt] (n7) at (0.6,0.5){};
\node[circle,fill,inner sep=1pt] (n8) at (1.4,0){};
\node[circle,fill,inner sep=1pt] (n9) at (-1.5,0){};
\draw [-] (n1) -- (n3);
\draw [-] (n1) -- (n7);
\draw [-] (n1) -- (n8);
\draw [-] (n3) -- (n7);
\draw [-] (n7) -- (n8);
\draw [-] (n4) -- (n2);
\draw [-] (n6) -- (n4);
\draw [-] (n4) -- (n9);
\draw [-] (n3) -- (n6);
\draw [-] (n4) -- (n7);
\draw [-] (n4) -- (n8);
\draw [-] (n2) -- (n8);
\draw [-] (n2) -- (n5);
\draw [-] (n5) -- (n9);
\draw [-] (n6) -- (n9);
\draw [-] (n1) -- (n3);
\draw [-] (n1) -- (n3);
\draw [-] (n9) -- (-2.75,0);
\draw [-] (n8) -- (2.75,0);
\draw [-] (n5) -- (-1.5,-2.3);
\draw [-] (n2) -- (1.5,-2.3);
\draw (-4.7,-0.8) [latex-latex,thick] (-4.7,1) node (yaxis){} |- (-3.2,-0.8) node (xaxis){};
\draw (-4.5,0.8) -- (-4,0.8) -- (-4,0.3) -- (-3.5,0.3) -- (-3.5,-0.3)-- (-3,-0.3){};
\draw (-4.7+1.5,-0.8-3) [latex-latex,thick] (-4.7+1.5,1-3) node (yaxis) {} |- (-3.2+1.5,-0.8-3) node (xaxis){};
\draw (-4.5+1.5,0.8-3) -- (-4+1.5,0.8-3) -- (-4+1.5,0.3-3) -- (-3.5+1.5,0.3-3) -- (-3.5+1.5,-0.3-3)-- (-3+1.5,-0.3-3);
\draw (3.4,-0.8) [latex-latex,thick] (3.4,1) node (yaxis) {} |- (4.9,-0.8) node (xaxis){};
\draw (3.5,0.8) -- (4,0.8) -- (4,0.3) -- (4.5,0.3) -- (4.5,-0.3)-- (5,-0.3);
\draw (3.4-1.5,-0.8-3) [latex-latex,thick] (3.4-1.5,1-3) node (yaxis) {} |- (4.9-1.5,-0.8-3) node (xaxis) {};
\draw (3.5-1.5,0.8-3) -- (4-1.5,0.8-3) -- (4-1.5,0.3-3) -- (4.5-1.5,0.3-3) -- (4.5-1.5,-0.3-3)-- (5-1.5,-0.3-3);
\end{tikzpicture}
\caption{Proposed TSO-DSO coordination approach}
\label{fig:paw_approach}
\end{figure}

In summary, we approximate the response of the distribution networks using a learning strategy that combines a K-nearest neighbor algorithm and a curve-fitting methodology. The proposed learning approach is simple and fast, as well as it does not suffer from high data requirements for training and offers explainability of the results, unlike  black-box approaches, e.g. based on deep learning. The operation of the transmission network is obtained by assuming that each distribution network reacts to prices according to the obtained step-wise non-increasing functions as illustrated in Fig. \ref{fig:paw_approach}. Mathematically, the operation of the transmission network is obtained by solving the following optimization problem: 
\begin{subequations}
\begin{align}
    & \max_{\Phi^T_t,p^N_{kt},p^B_{ktb}} \sum_{b \in \mathcal{B},k\in K ^T} \hspace{-2mm}\overline{u}^B_{ktb} p^B_{ktb} + \hspace{-2mm} \sum_{j\in J^T}u_{jt}(p^D_{jt}) - \hspace{-2mm} \sum_{i\in I^T} c_i(\pg) \label{eq:inverse_of}\\
    & \text{s.t.}\quad p^N_{kt} = \overline{p}^B_{kt0} + \sum_{b \in \mathcal{B}} p^B_{ktb}, \; \forall k \in K ^T \\
    & \phantom{s.t.}\quad 0 \leq p^B_{ktb} \leq \overline{p}^B_{ktb}-\overline{p}^B_{kt(b-1)}, \; \forall b \in \mathcal{B}, k\in K^T \\
    & \phantom{s.t.}\quad \eqref{eq:transmission_bal}- \eqref{eq:transmission_maxflow}
\end{align} \label{eq:inverse}
\end{subequations}

The proposed approach has several advantages. First, while the SB and PAG approaches disregard, respectively, the impact of network limits or economic signals on the response of distribution networks, the PAW approach is aware of both effects. Second, like the PAG approach, this method only requires historical LMPs  and power flows  at the substations and, therefore, detailed information about the distribution network parameters is not required. Third, the response of each distribution network to prices is modeled by a step-wise decreasing function that can be directly included in existing market-clearing mechanisms without  additional modifications. Besides, unlike other decomposition/decentralized approaches, the one we propose is not an iterative method and, therefore, is immune to convergence issues. Since the proposed method basically relies on a learning task, its performance highly depends on the quality of the input historical data. Hence, the dataset should be continuously updated to include the most recent operating conditions and exclude the oldest ones.

To conclude this section, Table \ref{tab:method_comparison} summarizes the main features of the four approaches discussed above. If compared with the benchmark, the three alternative approaches involve lower computational burdens through different approximation strategies. The next section describes the methodology to quantify the impact of such approximations on the optimal operation of the transmission electricity network.

\begin{table}[]
    \centering
    \caption{Qualitative comparison of TSO-DSOs coordination approaches}
    \begin{tabular}{lcccc}
    \toprule
    & BN & SB & PAG & PAW \\
    \midrule
    Network-aware & X & & X & X \\
    Price-aware & X & X & & X \\
    Historical data & & & X & X \\
    Seamless market integration & & & X & X \\
    Computational burden & High & Low & Low & Low \\
    \bottomrule
    \end{tabular}
    \label{tab:method_comparison}
\end{table}

\section{Evaluation procedure}\label{sec:CP}

While existing decomposition/decentralized approaches are able to yield the optimal coordination decisions, the proposed methodology may lead to suboptimal decisions caused by incompleteness or  inaccuracies of the learning task. In this section, we present the evaluation procedure to quantify the impact of these suboptimal decisions in terms of power imbalance and social welfare losses. We compare such measures with those obtained by the other three methods described in Section \ref{sec:methodology}. To that end, we proceed as follows:
\begin{enumerate}
    \item Solve problems \eqref{eq:benchmark}, \eqref{eq:dso2tso}, \eqref{eq:average} or \eqref{eq:inverse} using the modeling of the distribution networks derived from  the BN, SB, PAG or PAW approaches. LMPs at each substation $\lambda_{kt}$ are obtained as the dual variable of the balance equation \eqref{eq:transmission_bal}. The sum of the approximated consumption by all distribution networks is denoted as $\widehat{P}^N_t$.
    \item Model \eqref{eq:distribution} is solved for each distribution network $k$ after fixing LMPs at the substations to those obtained in Step 1). As such, we compute the actual response of the distribution networks considering all physical and economic information, denoted as $P^N_t$. Optimal values of objective function \eqref{eq:distribution_of} provide the social welfare achieved by each distribution network for the electricity prices computed in Step 1). We denote the sum of the social welfare of all distribution networks as $SW^D_t$.
    \item Quantify the power imbalance caused by the different distribution network approximations as $\Delta_t = 100 |\widehat{P}^N_t-P^N_t|/P^N_t$. Note that such power imbalances must be handled by flexible power resources able to adapt their generation or consumption in real-time.
    \item Model \eqref{eq:transmission} is solved by setting the electricity imported by each distribution network to the quantity obtained in Step 2). The output of this model represents the real-time re-dispatch of generating units connected to the transmission network to ensure the power system balance. The optimal value of \eqref{eq:transmission_of} provides the realized social welfare of the transmission network denoted as $SW^T_t$. We emphasize that this social welfare is computed as if all generating units and consumers at the transmission network could instantly adapt to any unexpected power imbalance coming from the distribution networks ($\Delta_t$) without any extra cost for the deployment of such unrealistic flexible resources. This means that we are underestimating the social welfare loss caused by these power imbalances.
    \item Compute the total realized social welfare of the power system as $SW_t = SW^D_t+SW^T_t$. 
\end{enumerate}

\section{Simulation results}\label{sec:SR}

We consider the 118-bus, 186-line transmission network from \cite{Pena2018}. Each  transmission-level load is replaced with a 32-bus radial distribution network, which hosts eight solar generating units, see data in \cite{Hassan2018b, Hassan2018}. That is, the power system includes 3030 buses ($118+91\times32$),  3098 lines ($186+91\times32$), thermal and wind power plants  connected to 43 transmission buses, solar generating  units  connected to 728 distribution buses ($91\times8$), and electricity consumers located at 2912 distribution buses ($91\times32$). Each consumer is assumed to react to the electricity price as depicted in Fig. \ref{fig:flex}. The installed capacity of thermal, solar and wind generating units is 17.3GW, 2.5GW and 2.5GW, respectively, while the peak demand is 18GW. Finally, time-varying capacity factors for all consumers, wind and solar generation in the same distribution network are assumed equal. While all distribution networks have the same topology and the same location of loads and solar power generating units, we scale their total demand from 12MW to 823MW to match the transmission demand given in \cite{Pena2018}. We also scale the original values of branch resistances and reactances inversely proportional to the  peak demand within each distribution network. All data used in this case study is available in \cite{118TN33DN}. Simulations have been run on a Linux-based server with one CPU clocking at 2.6 GHz and 2 GB of RAM using CPLEX 12.6 under Pyomo~5.2.

As discussed in Section \ref{sec:methodology}, the analyzed methods  differ in their ability to account for the impact of physical limits and economic signals on the response of active distribution networks. For instance, if distribution voltage limits never become activated, then the SB approach would provide results quite close to those of the benchmark approach BN. Conversely, if distribution voltages reach their security limits, the PAG and PAW methods are expected to outperform SB. 
In order to investigate the impact of voltage congestion on the  performance of each approach, we vary the resistances and reactances of branches of the distribution networks as indicated in \eqref{eq:change_rx}, where $r^0,x^0$ are the base-case values provided in \cite{118TN33DN}, and parameter $\eta$ is changed from 0.67 to 1.33, i.e., a 33\% lower and greater than the initial values:
\begin{equation}
r = \eta r^0 \qquad \qquad x = \eta x^0 \label{eq:change_rx}
\end{equation}
Additionally, we use parameter $\delta$ to model each flexible consumer, which is randomly generated for the 2912 loads following a uniform probability distribution in $[0.5-0.75]$. Besides, we set $\overline{\lambda}=25$, $\underline{\lambda}=10$, $\overline{v}_n=1.05$, and $\underline{v}_n=0.95$.

\begin{figure}
\centering
\includegraphics[scale=0.55]{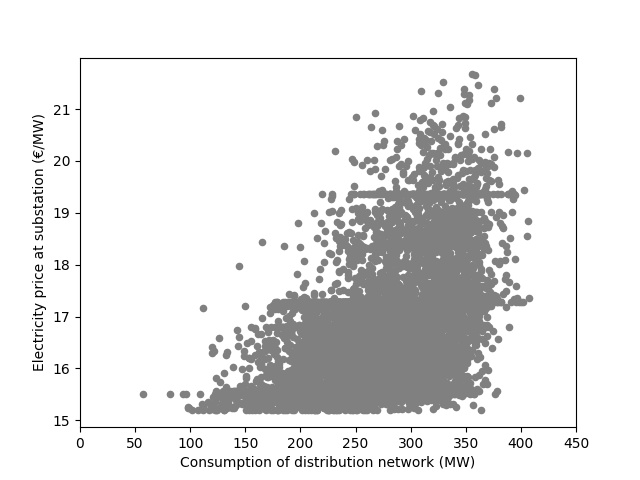} \caption{Price response of one distribution network for $\eta=1$}
\label{fig:response}
\end{figure}

The PAG and PAW approaches require access to historical data. In this case study,  historical data is  generated by solving the  BN model \eqref{eq:benchmark} for 8760 hours of a given year. Each hour is characterized by different baseline demands of flexible consumers along with the wind and solar capacity factors throughout the system. Values for wind, solar and baseline demands are taken from \cite{Pena2018} and are available in \cite{118TN33DN}. For illustration, Fig. \ref{fig:response} plots the price response of one distribution network for $\eta=1$ including the 8760 time periods. For simplicity, changes in the topology of the transmission and distribution networks are disregarded. The learning-based PAG and PAW approaches use the demand and renewable capacity factors at each distribution network as contextual information to learn its response. Also, the number of neighbors for the $K$-NN learning methodology is set to 100. Finally, the maximum number of blocks for the  bidding curves learned by the PAW approach is equal to ten. For the sake of comparison, each of the four approaches uses the same test set that includes 100 randomly selected hours of the year. 

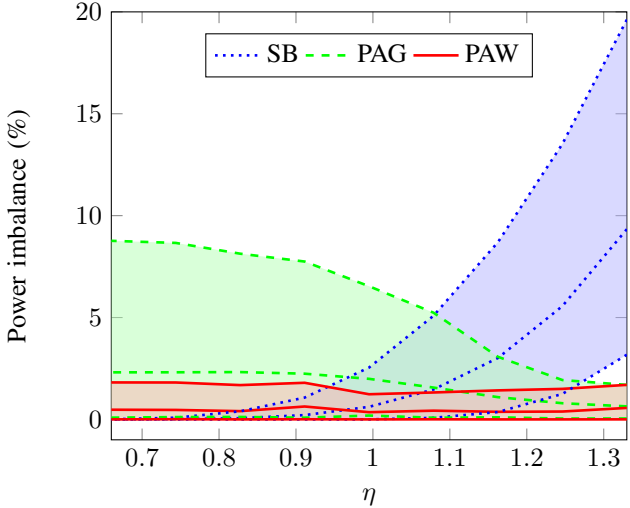
\begin{figure}
\centering
\begin{tikzpicture}[scale=1]
	\begin{axis}[	
    xmin = 0.66,
    xmax = 1.33,
    ymin = -1,
    ymax = 20,
    legend style={at={(0.5,0.95)},anchor=north,legend cell align=left,legend columns=4},	
	clip marker paths=true,	
	xlabel = $\eta$,
	ylabel = Power imbalance (\%)
    ]	
  \addplot[line width=1pt,draw=blue,dotted] table [x=eta, y=demrelerr, col sep=comma,discard if not={mode}{dso2tso}]{resultsMF.csv}; \addlegendentry{SB}
  \addplot[line width=1pt,draw=green,dashed] table [x=eta, y=demrelerr, col sep=comma,discard if not={mode}{average}]{resultsMF.csv}; \addlegendentry{PAG}
  \addplot[line width=1pt,draw=red] table [x=eta, y=demrelerr, col sep=comma,discard if not={mode}{inverse}]{resultsMF.csv}; \addlegendentry{PAW}
  \addplot[line width=1pt,draw=blue,name path=sb95,dotted] table [x=eta, y=dem95, col sep=comma,discard if not={mode}{dso2tso}]{resultsMF.csv}; 
  \addplot[line width=1pt,draw=green,name path=pag95,dashed] table [x=eta, y=dem95, col sep=comma,discard if not={mode}{average}]{resultsMF.csv}; 
  \addplot[line width=1pt,draw=red,name path=paw95] table [x=eta, y=dem95, col sep=comma,discard if not={mode}{inverse}]{resultsMF.csv}; 
  \addplot[line width=1pt,draw=blue,name path=sb05,dotted] table [x=eta, y=dem05, col sep=comma,discard if not={mode}{dso2tso}]{resultsMF.csv}; 
  \addplot[line width=1pt,draw=green,name path=pag05,dashed] table [x=eta, y=dem05, col sep=comma,discard if not={mode}{average}]{resultsMF.csv}; 
  \addplot[line width=1pt,draw=red,name path=paw05] table [x=eta, y=dem05, col sep=comma,discard if not={mode}{inverse}]{resultsMF.csv}; 
  \addplot[fill=blue!30,fill opacity=0.5] fill between[of=sb05 and sb95];
  \addplot[green!30,fill opacity=0.5] fill between[of=pag05 and pag95];
  \addplot[red!30,fill opacity=0.5] fill between[of=paw05 and paw95];
];
\end{axis}	
\end{tikzpicture} \caption{Impact of distribution network congestion on power imbalance.} \label{fig:imbalance}
\end{figure}

Using the results of these 100 hours, Fig. \ref{fig:imbalance} plots, for each approach, a shaded area ranging from the 5\% to the 95\% percentile of the relative power imbalance $\Delta_t$ as a function of parameter $\eta$. The average of the power imbalance is also displayed. Naturally, due to its completeness, the benchmark method does not incur any power imbalance and therefore, the results delivered by this method are not included in Fig. \ref{fig:imbalance}. Low values of $\eta$ reduce voltage congestion  at the distribution networks and, therefore, their response is mainly driven by electricity prices at the substations. In such cases, the SB approach outperforms the PAG approach and yields power imbalances close to 0\%. For small values of $\eta$, the proposed PAW approach yields higher power imbalances than SB. However, this difference could be narrowed by approximating the response of the distribution networks with more than ten blocks. Conversely, high values of $\eta$ translates into congested distribution networks in which the dispatch of small consumers and distributed generators is heavily constrained by technical limits. In these circumstances, electricity prices at the substations have a reduced impact on the the response of the distribution network and then, the power imbalance of the SB approach is significantly greater than that of the PAG approach. Quantitatively, the proposed methodology PAW achieves average power imbalances below 0.7\% for any value of $\eta$.

When comparing SB, PAG and PAW, we shoud also keep in mind that their integration into current market-clearing mechanisms are not comparable. Implementing the SB approach would require modifying existing market rules so that distributed generators and small consumers could directly submit their offers and bids. On the other hand, the PAG and PAW comply with these rules since active distribution networks are modeled as fix loads or in the form of step-wise bidding curves, respectively. 

\begin{figure}
\centering
\begin{tikzpicture}[scale=1]
	\begin{axis}[	
    xmin = 0.66,
    xmax = 1.33,
    ymin = 3,
    ymax = 12,
    legend style={at={(0.5,0.95)},anchor=north,legend cell align=left,legend columns=4},	
	clip marker paths=true,	
	xlabel = $\eta$,
	ylabel = Demand variation wrt baseline (\%)
    ]	
  \addplot[line width=1pt,draw=blue] table [x=eta, y=realflex, col sep=comma,discard if not={mode}{dso2tso}]{resultsMF.csv}; \addlegendentry{SB}
  \addplot[line width=1pt,draw=green] table [x=eta, y=realflex, col sep=comma,discard if not={mode}{average}]{resultsMF.csv}; \addlegendentry{PAG}
  \addplot[line width=1pt,draw=red] table [x=eta, y=realflex, col sep=comma,discard if not={mode}{inverse}]{resultsMF.csv}; \addlegendentry{PAW}
  \addplot[line width=1pt,draw=black] table [x=eta, y=realflex, col sep=comma,discard if not={mode}{normal}]{resultsMF.csv}; \addlegendentry{BN}
  \addplot[line width=1pt,draw=blue,dashed] table [x=eta, y=foreflex, col sep=comma,discard if not={mode}{dso2tso}]{resultsMF.csv}; 
  \addplot[line width=1pt,draw=green,dashed] table [x=eta, y=foreflex, col sep=comma,discard if not={mode}{average}]{resultsMF.csv}; 
  \addplot[line width=1pt,draw=red,dashed] table [x=eta, y=foreflex, col sep=comma,discard if not={mode}{inverse}]{resultsMF.csv}; 
  \addplot[line width=1pt,draw=black,dashed] table [x=eta, y=foreflex, col sep=comma,discard if not={mode}{normal}]{resultsMF.csv}; 
];
\end{axis}	
\end{tikzpicture} \caption{Estimated (dashed) and observed (bold) flexibility of active distribution networks.} \label{fig:usedflex}
\end{figure}
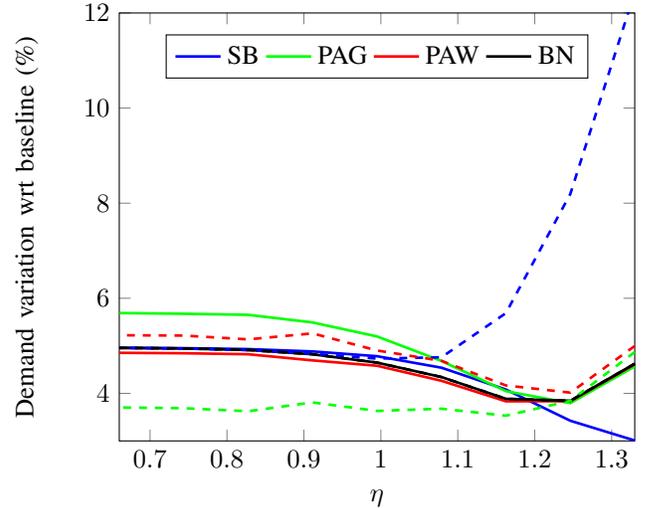

Power imbalances of Fig. \ref{fig:imbalance} are also explained by the incorrect estimation of the flexibility provided by the distribution networks made by the different approaches. To illustrate this effect, we compute the relative difference between the approximate consumption of distribution networks $\widehat{P}^N_t$ and the baseline consumption, which is plotted in dashed lines in Fig. \ref{fig:usedflex}. The bold lines represent the relative difference between the actual consumption of distribution networks $P^N_t$  and the same baseline demand. First, it can be observed that flexible customers allow for aggregate demand variations that range from 4\% to 12\%, on average. It can also be observed that PAG underestimates the flexibility provided by distribution networks for low congestion levels. Conversely, the SB approach overestimates the available flexibility of distributed networks when their operation is mainly driven by physical constraints. Finally, the proposed PAW approach is able to operate the transmission network with a very realistic estimation of the flexibility of the distribution networks, which is, in turn, very close to the actual flexibility levels determined by the centralized benchmark approach.

\begin{figure}
\centering
\begin{tikzpicture}[scale=1]
	\begin{axis}[	
    xmin = 0.66,
    xmax = 1.33,
    ymin = -0.2,
    ymax = 5,
    legend style={at={(0.5,0.95)},anchor=north,legend cell align=left,legend columns=4},	
	clip marker paths=true,	
	xlabel = $\eta$,
	ylabel = Social welfare loss (\%)
    ]	
  \addplot[line width=1pt,draw=blue,dotted] table [x=eta, y=swrelerr, col sep=comma,discard if not={mode}{dso2tso}]{resultsMF.csv}; \addlegendentry{SB}
  \addplot[line width=1pt,draw=green,dashed] table [x=eta, y=swrelerr, col sep=comma,discard if not={mode}{average}]{resultsMF.csv}; \addlegendentry{PAG}
  \addplot[line width=1pt,draw=red] table [x=eta, y=swrelerr, col sep=comma,discard if not={mode}{inverse}]{resultsMF.csv}; \addlegendentry{PAW}
  \addplot[line width=1pt,draw=blue,name path=sb95,dotted] table [x=eta, y=sw95, col sep=comma,discard if not={mode}{dso2tso}]{resultsMF.csv}; 
  \addplot[line width=1pt,draw=green,name path=pag95,dashed] table [x=eta, y=sw95, col sep=comma,discard if not={mode}{average}]{resultsMF.csv}; 
  \addplot[line width=1pt,draw=red,name path=paw95] table [x=eta, y=sw95, col sep=comma,discard if not={mode}{inverse}]{resultsMF.csv}; 
  \addplot[line width=1pt,draw=blue,name path=sb05,dotted] table [x=eta, y=sw05, col sep=comma,discard if not={mode}{dso2tso}]{resultsMF.csv}; 
  \addplot[line width=1pt,draw=green,name path=pag05,dashed] table [x=eta, y=sw05, col sep=comma,discard if not={mode}{average}]{resultsMF.csv}; 
  \addplot[line width=1pt,draw=red,name path=paw05] table [x=eta, y=sw05, col sep=comma,discard if not={mode}{inverse}]{resultsMF.csv}; 
  \addplot[fill=blue!30,fill opacity=0.5] fill between[of=sb05 and sb95];
  \addplot[green!30,fill opacity=0.5] fill between[of=pag05 and pag95];
  \addplot[red!30,fill opacity=0.5] fill between[of=paw05 and paw95];
];
\end{axis}	
\end{tikzpicture} \caption{Impact of distribution network congestion on social welfare.} \label{fig:socialwelf}
\end{figure}
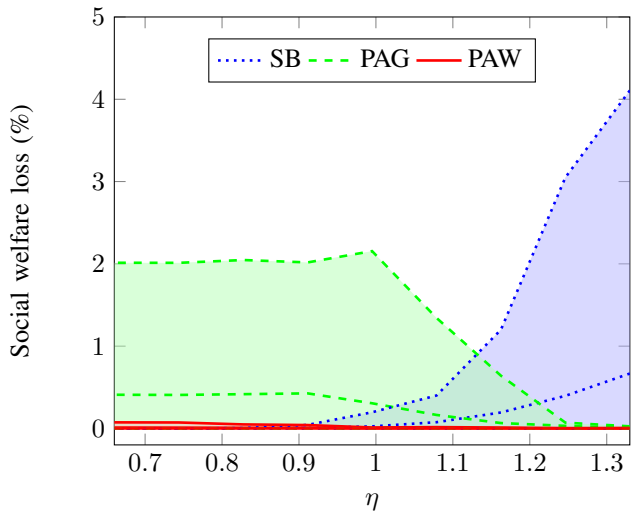

Similarly to Fig. \ref{fig:imbalance}, Fig. \ref{fig:socialwelf} plots the mean and the 5\% and 95\% percentiles of the social welfare loss with respect to the BN approach. Aligned with power imbalance results, the social welfare losses under the SB and PAG approaches are linked to high and low values of parameter $\eta$, respectively. More importantly, while the social welfare loss may reach values of 2\% and 4\% for the SB and PAG approaches, in that order, for some of the 100 hours analyzed, the PAW approach keeps this value below 0.1\% for any network congestion level. That is, the proposed methodology to integrate transmission and distribution networks achieves the same social welfare as the BN for a wide range of power system conditions (described by the different demand and renewable capacity factors of the 100 hours) and network congestion of the distribution systems (modeled by parameter $\eta$). For completeness, Table~\ref{tab:socialwelfares} provides the sum of the social welfare for the 100 hours of the test set for some values of parameter $\eta$.

It is also important to remark that social welfare losses in Fig. \ref{fig:socialwelf} are computed assuming that all generating units and consumers at the transmission network can react instantaneously to any real-time power imbalance without involving extra regulations costs. Therefore, these results are a lower bound of the actual social welfare losses that would happen in a more realistic setup in which flexibility resources are both limited and expensive.

\begin{table}[]
    \centering
    \caption{Social welfare results in k\euro}
\begin{tabular}{ccccc}
\toprule
  $\eta$ &  BN & SB & PAG & PAW \\
\midrule
 0.66 &  3226.8 &  3226.8 &  3216.8 &  3226.3 \\
 1.00 &  3214.5 &  3213.5 &  3207.5 &  3214.4 \\
 1.33 &  2828.8 &  2814.1 &  2828.7 &  2828.7 \\
\bottomrule
\end{tabular}
\label{tab:socialwelfares}
\end{table}

\begin{table}[]
    \centering
    \caption{Allocation of average social welfare loss (in percent with respect to BN)  between transmission and distribution}
    \begin{tabular}{ccccccc}
    \toprule
    & \multicolumn{2}{c}{SB} & \multicolumn{2}{c}{PAG} & \multicolumn{2}{c}{PAW} \\
    \midrule
    $\eta$& TSO & DSO & TSO & DSO & TSO & DSO \\
    \midrule
    0.66 & 0.00\% & 0.00\% & -4.98\% & 5.39\% & 0.46\% & -0.45\% \\
    1.00 & -1.38\% & 1.41\% & -4.17\% & 4.47\% & 0.04\% & -0.04\%\\
    1.33 & -12.00\% & 12.67\% & -0.54\% & 0.55\% & -0.10\% & 0.10\%\\
    \bottomrule
    \end{tabular}
    \label{tab:socialwelfloss_alloc}
\end{table}

Table \ref{tab:socialwelfloss_alloc} shows how the average relative social welfare loss (as illustrated in Fig.~\ref{fig:socialwelf}) is apportioned between the transmission and distribution systems, for various congestion levels $\eta$. Notably, the average loss in the SB and PAG cases  disproportionally affects the transmission and distribution networks. Actually, there is a substantial net transfer of welfare from DSOs to the TSO. That is, the SB and PAG approaches delegate the bulk of the costs of dealing with distribution congestion to the distributed energy resources themselves, which certainly puts into question the ability of these methods to effectively integrate distribution into transmission operations. In contrast, the proposed PAW approach considerably mitigates this effect, or even reverses it, thus ensuring that distribution issues are also taken care of by transmission resources.  

Looking at Figures \ref{fig:imbalance}, \ref{fig:usedflex} and \ref{fig:socialwelf}, we can conclude that the effectiveness of the proposed PAW method with respect to other approaches depends on the network characteristics, and more particularly, on the congestion level of the distribution networks. Indeed, if the distribution networks never experience voltage congestion, SB outperforms PAW in terms of power imbalance and social welfare loss. On the other hand, for highly congested networks, PAG and PAW provide almost identical results. In conclusion, the use of the proposed PAW approach is more relevant for those systems in which the level of congestion of distribution networks vary significantly depending on the operating conditions.

Finally, Table \ref{tab:computational_times} compares the maximum, average and minimum computational times for the four approaches. The average speedup factor between each method and the benchmark is also provided in the last column. Due to the high number of variables and constraints of model \eqref{eq:dso2tso}, the SB's speedup factor is relatively low. In contrast, since PAG and PAW characterize the response of each distribution network through a constant value or a step-wise bidding curve, respectively, the computational savings are more substantial.

\begin{table}[]
    \centering
    \caption{Computational time results}
\begin{tabular}{lcccc}
\toprule
& Min time (s) & Average time (s) & Max time (s) & Speedup \\
\midrule
BN & 0.88 & 1.66 & 15.27 &  - \\
SB & 0.17 & 0.30 &  1.83 &  5.5x \\
PAG & 0.02 & 0.03 &  0.20 & 55.3x \\
PAW & 0.04 & 0.06 &  0.38 & 27.7x \\
\bottomrule
\end{tabular}
\label{tab:computational_times}
\end{table}

\section{Conclusion}\label{sec:conc}

Motivated by the proliferation of distributed energy resources, new TSO-DSO coordination strategies are required to take full advantage of these resources in the operation of the transmission system. Existing decomposition/decentralized methods are able to yield the same operating decisions as centralized benchmarks at lower computational costs. However, these approaches require access to proprietary information or fail-safe real-time communication infrastructures. Alternatively, our approach only uses offline historical data at substations to learn the price response of the distribution networks in the form of a non-increasing bidding curve that can be easily embedded into current procedures for transmission operations. In addition, this data set can be enriched with some covariates that have predictive power on the response of the distribution networks.

We have benchmarked our approach against an idealistic model that fully centralizes the coordination of distribution and transmission operations. We have also compared it with other approximate approaches that either ignore the technical constraints of the distribution networks or the price-sensitivity of DERs. The conducted numerical experiments reveal that our approach systematically delivers small differences with respect to the fully centralized benchmark in terms of power imbalances and social welfare regardless of the level of congestion of the distribution grids. In return, our approach is computationally affordable and consistent with current market practices, and allows for decentralization.

Future work should be directed to assessing whether these results remain valid, and to which extent, for meshed distribution networks and DERs with more complex price responses, e.g. thermostatically controlled loads. Furthermore, in this research, we have only considered contextual information pertaining to continuous random variables (electricity demand, renewable power generation, etc.). Therefore, a relevant avenue for future research is to extend the proposed approach to work with binary variables too, such as those describing the network topology or the in-service/out-of-service status of some network components.

\bibliographystyle{IEEEtran}
\bibliography{references}

\end{document}